\documentclass[preprint]{elsarticle}

\usepackage{lineno,hyperref}
\modulolinenumbers[5]

\journal{Journal of \LaTeX\ Templates}

\usepackage{amsmath,amssymb}
\usepackage{dsfont}
\usepackage[shortlabels]{enumitem}
\newtheorem{thm}{Theorem}[section]
\newtheorem{defn}[thm]{Definition}
\newtheorem{lem}[thm]{Lemma}

\newtheorem{rem}[thm]{Remark}
\newtheorem{exa}[thm]{Example}
\newtheorem{ass}[thm]{Assumption}

\allowdisplaybreaks

\usepackage{color,soul}
\usepackage{mathtools}

\usepackage{graphicx}

\newcounter{enum}

\PassOptionsToPackage{unicode}{hyperref}
\PassOptionsToPackage{naturalnames}{hyperref}

\begin{document}

\begin{frontmatter}

\title{Bayesian inverse problems with partial observations}

\author[mymainaddress]{Shota Gugushvili\corref{mycorrespondingauthor}}
\cortext[mycorrespondingauthor]{Corresponding author}
\ead{shota.gugushvili@math.leidenuniv.nl}

\author[mymainaddress]{Aad W. van der Vaart}
\ead{avdvaart@math.leidenuniv.nl}

\author[mymainaddress]{Dong Yan}
\ead{d.yan@math.leidenuniv.nl}

\address[mymainaddress]{Mathematical Institute, Faculty of Science, Leiden University, P.O. Box 9512, 2300 RA Leiden, The Netherlands}

\begin{abstract}
We study a nonparametric Bayesian approach to linear inverse problems under discrete observations. We use the discrete Fourier transform to convert our model into a truncated Gaussian sequence model, that is closely related to the classical Gaussian sequence model. Upon placing the truncated series prior on the unknown parameter, we show that as the number of observations $n\rightarrow\infty,$ the corresponding posterior distribution contracts around the true parameter at a rate depending on the smoothness of the true parameter and the prior, and the ill-posedness degree of the problem. Correct combinations of these values lead to optimal posterior contraction rates (up to logarithmic factors). Similarly, the frequentist coverage of Bayesian credible sets is shown to be dependent on a combination of smoothness of the true parameter and the prior, and the ill-posedness of the problem. Oversmoothing priors lead to zero coverage, while undersmoothing priors produce highly conservative results. Finally, we illustrate our theoretical results by numerical examples.
\end{abstract}


\end{frontmatter}

\section{Introduction}
\label{sec:Introduction}

Linear inverse problems have been studied since long in the statistical and numerical analysis literature; see, e.g., \cite{alquier2011inverse}, \cite{bissantz2007inverseregularization}, \cite{cavalier2008nonparametricinverseproblems}, \cite{cavalier2002sharp}, \cite{cohen2004adaptive_galerkin}, \cite{donoho1995WVD}, \cite{kaipio2006statistical}, \cite{kirsch2011introduction}, \cite{wahba1977integraloperator}, and references therein. Emphasis in these works has been on the {signal-in-white noise model},
\begin{equation}
\label{eq:general_inverse_problem}
Y = Af + \varepsilon W,
\end{equation}
where the parameter of interest $f$ lies in some infinite-dimensional function space, $A$ is a linear operator with values in a possibly different space, $W$ is white noise, and $\varepsilon$ is the noise level. Applications of linear inverse problems include, e.g., computerized tomography, see \cite{natterer2001mathematics_CT}, partial differential equations, see \cite{isakov2013inverse}, and scattering theory, see \cite{colton2012inverse_acoustic_electromagnetic}.

Arguably, in practice one does not have access to a full record of observations on the unknown function $f$ as in the idealised model \eqref{eq:general_inverse_problem}, but rather one indirectly observes it at a finite number of points. This statistical setting can be conveniently formalised as follows: let the signal of interest $f$ be an element in a Hilbert space $H_1$ of functions defined on a compact interval $[0,1]$. The {forward operator} $A$ maps $f$ to another Hilbert space $H_2$. We assume that $H_1$, $H_2$ are subspaces of $L^2([0,1])$, typically collections of functions of certain smoothness as specified in the later sections, and that the design points are chosen deterministically,
\begin{equation}
\label{eq:design_points}
\left\{x_i = \frac{i}{n}\right\}_{i=1,\cdots,n}.
\end{equation}
Assuming continuity of $Af$ and defining
\begin{equation}
\label{eq:general_inverse_problem_discrete}
Y_i = Af(x_i) + \xi_i, \quad i = 1,\cdots, n,
\end{equation}
with $\xi_i$ i.i.d.\ standard Gaussian random variables, our observations are the pairs $(x_i, Y_i)_{i\leq n},$ and we are interested in estimating $f$. A prototype example we think of is the case when $A$ is the solution operator in the Dirichlet problem for the heat equation acting on the initial condition $f$; see Example \ref{exa:heat_equation} below for details.

Model \eqref{eq:general_inverse_problem_discrete} is related to the inverse regression model studied e.g. in \cite{birke2010inverse_regression} and \cite{bissantz2012inverse_regression}. Although the setting we consider is somewhat special, our contribution is arguably the first one to study from a theoretical point of view a nonparametric Bayesian approach to estimation of $f$ in the inverse problem setting with partial observations (see \cite{aad:book:2017} for a monographic treatment of modern Bayesian nonparametrics).  In the context of the signal-in-white noise model \eqref{eq:general_inverse_problem}, a nonparametric Bayesian approach has been studied thoroughly in \cite{knapik2011bayesianmild} and \cite{knapik2013bayesianextreme}, and techniques from these works will turn out to be useful in our context as well. Our results will deal with derivation of posterior contraction rates and study of asymptotic frequentist coverage of Bayesian credible sets. A posterior contraction rate can be thought of as a Bayesian analogue of a convergence rate of a frequentist estimator, cf.~\cite{vandervaart2000posterior} and \cite{aad:book:2017}. Specifically, we will show that as the sample size $n\rightarrow\infty$, the posterior distribution concentrates around the `true' parameter value, under which data have been generated, and hence our Bayesian approach is consistent and asymptotically recovers the unknown `true' $f$. The rate at which this occurs will depend on the smoothness of the true parameter and the prior and the ill-posedness degree of the problem. Correct combinations of these values lead to optimal posterior contraction rates (up to logarithmic factors). Furthermore, a Bayesian approach automatically provides uncertainty quantification in parameter estimation through the spread of the posterior distribution, specifically by means of posterior credible sets. We will give an asymptotic frequentist interpretation of these sets in our context. In particular, we will see that the frequentist coverage will depend on a combination of smoothness of the true parameter and the prior, and the ill-posedness of the problem. Oversmoothing priors lead to zero coverage, while undersmoothing priors produce highly conservative results.

The article is organized as follows:
in Section~\ref{sec:Problem}, we give a detailed description of the problem, introduce the singular value decomposition and convert the model \eqref{eq:general_inverse_problem_discrete} into an equivalent truncated sequence model that is better amenable to our theoretical analysis. We show how a Gaussian prior in this sequence model leads to a Gaussian posterior and give an explicit characterisation of the latter.
Our main results on posterior contraction rates and Bayesian credible sets are given in Section~\ref{sec:MainResults}, followed by simulation examples  in Section~\ref{sec:SimulationExamples} that illustrate our theoretical results. Section~\ref{sec:proofs} contains the proofs of the main theorems, while the technical lemmas used in the proofs are collected in Section~\ref{sec:appendix}.

\subsection{Notation}
The notational conventions we use in this work are the following: definitions are marked by the $:=$ symbol; $|\cdot |$ denotes the absolute value and $\|\cdot\|_{H}$ indicates the norm related to the space $H$; $\langle \cdot,\cdot\rangle_{H}$ is understood as the canonical inner product in the inner product space $H$; subscripts are omitted when there is no danger of confusion; $\mathcal{N}(\mu,\Sigma)$ denotes the Gaussian distribution with mean $\mu$ and covariance operator $\Sigma$; subscripts $\mathcal{N}_n$ and $\mathcal{N}_H$ may be used to emphasize the fact that the distribution is defined on the space $\mathbb{R}^n$ or on the abstract space $H$; $\operatorname{Cov}(\cdot,\cdot)$ denotes the covariance or the covariance operator, depending on the context; for positive sequences $\{a_n\},$ $\{b_n\}$ of real numbers, the notation $a_n\lesssim b_n$ and $a_n \gtrsim b_n$ mean respectively that there exist positive constants $C_1,C_2$ independent of $n,$ such that $a_n\leq C_1b_n$ or $a_n \geq C_2 b_n$ hold for all $n$; finally, $a_n\asymp b_n$ indicates that the ratio $a_n/b_n$ is asymptotically bounded from zero and infinity, while $a_n \sim b_n$ means $a_n/b_n \to 1$ as $n \to \infty$.

\section{Sequence model}
\label{sec:Problem}

\subsection{Singular value decomposition}
\label{subsec:SVD}

We impose a common assumption on the forward operator $A$ from the literature on inverse problems, see, e.g., \cite{alquier2011inverse}, \cite{bissantz2007inverseregularization} and \cite{cavalier2008nonparametricinverseproblems}.

\begin{ass}
\label{ass:OperatorA}
Operator $A$ is injective and compact. 
\end{ass}

It follows that $A^*A$ is also compact and in addition self-adjoint. Hence, by the spectral theorem for self-adjoint compact operators, see \cite{conway1990courseonFA}, we have a representation
$
A^*Af = \sum_{k\in\mathbb{N}} a_k^2 f_k \varphi_k,
$
where $\{\varphi_k\}$ and $\{a_k\}$ are the {eigenbasis} on $H_1$ and {eigenvalues}, respectively, (corresponding to the operator $A^*A$), and $f_k = \langle f,\varphi_k\rangle$ are the {Fourier coefficients} of $f$. This decomposition of $A^*A$ is known as the {singular value decomposition (SVD)}, and $\{a_k\}$ are also called {singular values}.

It is easy to show that the {conjugate basis} $\psi_k:=A\varphi_k/a_k$  of the orthonormal basis $\{\varphi_k\}_k$ is again an orthonormal system in $H_2$ and gives a convenient basis for $Range(A)$, the range of $A$ in $H_2$. Furthermore, the following relations hold (see \cite{alquier2011inverse}),
\begin{equation}
\label{eq:BasisMappingRelation}
A\varphi_k = a_k\psi_k,\ A^*\psi_k = a_k\varphi_k.
\end{equation}

Recall a standard result (see, e.g., \cite{haase2014functional}): a Hilbert space $H$ is isometric to $\ell^2,$ and Parseval's identity $\|f\|_{\ell^2}^2 :=\sum_k |f_k|^2 = \|f\|_{H}^2$ holds; here $f_k$ are the Fourier coefficients with respect to some known and fixed orthonormal basis.

We will employ the eigenbasis $\{\varphi_k\}$ of $A^*A$ to define the Sobolev space of functions. This will define the space in which the unknown function $f$ resides.

\begin{defn}
\label{def:SobolevSpace}
We say $f $ is in the Sobolev space $S^\beta$ with smoothness parameter $\beta\geq 0$, if
it can be written as
$
f = \sum_{k=1}^{\infty}f_{k}\varphi_k
$
{with}
$
f_k = \langle f,\varphi_k\rangle,
$
and if its norm
$
\|f\|_\beta:= \left(\sum_{k=1}^{\infty}f_k^2 k^{2\beta}\right)^{1/2}
$
is finite.
\end{defn}
\begin{rem}
The above definition agrees with the classical definition of the Sobolev space if the eigenbasis is the trigonometric basis, see, e.g., \cite{tsybakov2008introduction}. With a fixed basis, which is always the case in this article, 
one can identify the function $f$ and its Fourier coefficients $\{f_k\}$. Thus, we use $S^\beta$ to denote both the function space and the sequence space. For example, it is easy to verify that $S^0 = \ell^2$ (correspondingly $S^0 = L^2$), $S^\beta\subset \ell^2$ for any nonnegative $\beta$, and $S^\beta \subset \ell^1$ when $\beta > 1/2$.
\end{rem}

Recall that $A f = \sum a_i f_i \psi_i.$ Then we have $Af \in S^{\beta+p}$ if $a_k \asymp k^{-p},$ and $Af \in S^{\infty}:= \cap_{k \in \mathbb{N}} S^k,$ if $a_k$ decays exponentially fast.
Such a lifting property is beneficial in the forward problem, since it helps to obtain a smooth solution. However, in the context of inverse problems it leads to a difficulty in recovery of the original signal $f$, since information on it is washed out by smoothing. Hence, in the case of inverse problems one does not talk of the lifting property, but of ill-posedness, see \cite{cavalier2008nonparametricinverseproblems}.

\begin{defn}
\label{def:illposedness}
An inverse problem is called mildly ill-posed, if $a_k \asymp k^{-p}$  as $k \to \infty$, and extremely ill-posed, if $a_k \asymp e^{-k^s p}$ with $s \geq 1$ as $k \to \infty$, where $p$ is strictly positive in both cases.
\end{defn}

In the rest of the article, we will confine ourselves to the following setting.

\begin{ass}
The unknown true signal $f$ in \eqref{eq:general_inverse_problem_discrete} satisfies $f\in S^{\beta} \subset H_1$ for $\beta>0.$ Furthermore, the ill-posedness is of one of the two types in Definition~\ref{def:illposedness}.
\end{ass}
\begin{rem}
As an immediate consequence of the lifting property, we have $H_2 \subset H_1$.
\end{rem}

We conclude this section with two canonical examples of the operator $A$.

\begin{exa}[mildly ill-posed case: Volterra operator \cite{knapik2011bayesianmild}]
\label{exa:volterra_operator}
The classical Volterra operator $A:L^2[0,1] \to L^2[0,1]$ and its adjoint $A^*$ are
\[
Af(x) = \int_0^x f(s) \, ds,\quad A^* f(x) = \int_x^1 f(s) \, ds.
\]
The eigenvalues, eigenfunctions of $A^*A$ and the conjugate basis are given by 
\begin{align*}
a_i^2 =&\frac{1}{(i-1/2)^2 \pi^2},\\
\varphi_i(x) =& \sqrt{2} \cos ((i-1/2)\pi x),\\
\psi_i(x) =& \sqrt{2} \sin ((i-1/2)\pi x),
\end{align*}
for $i\geq 1$.
\end{exa}

\begin{exa}[extremely ill-posed case: heat equation \cite{knapik2013bayesianextreme}]
\label{exa:heat_equation}
Consider the Dirichlet problem for the heat equation:
\begin{equation}
\label{eq:heat_equation}
\begin{split}
\frac{\partial}{ \partial t} u(x,t) &= \frac{\partial^2}{ \partial x^2} u(x,t),\quad
u(x,0) = f(x), \\
u(0,t) &= u(1,t) = 0,\quad t\in[0,T],
\end{split}
\end{equation}
where $u(x,t)$ is defined on $[0,1]\times [0,T]$ and $f(x)\in L^2[0,1]$ satisfies $f(0) = f(1) = 0$.
The solution of \eqref{eq:heat_equation} is given by
\[
u(x,t) =\sqrt{2} \sum_{k=1}^{\infty} f_k e^{-k^2\pi^2 t} \sin(k\pi x) =: Af(x),
\]
where $\{f_k\}$ are the coordinates of $f$ in the basis $\{\sqrt{2} \sin(k\pi x)\}_{k\geq 1}$.

For the solution map $A$, the eigenvalues of $A^*A$ are $e^{-k^2\pi^2 t}$, the eigenbasis and conjugate basis coincide and $\varphi_k(x) =\psi_k(x) = \sqrt{2} \sin(k\pi x)$.
\end{exa}

\subsection{Equivalent formulation}
\label{subsec:gaussianInverse_equivalentFormulation}
In this subsection we develop a sequence formulation of the model \eqref{eq:general_inverse_problem_discrete}, which is very suitable for asymptotic Bayesian analysis. First, we briefly discuss the relevant results that provide motivation for our reformulation of the problem.

In Examples \ref{exa:volterra_operator} and \ref{exa:heat_equation}, the sine and cosine bases form the eigenbasis. In fact, the Fourier basis (trigonometric polynomials) frequently arises as an eigenbasis for various operators, e.g.\ in the case of differentiation, see \cite{efromovich1998differentiation}, or circular deconvolution, see \cite{cavalier2002sharp}. For simplicity, we will use Fourier basis as a primary example in the rest of the article. Possible generalization to other bases is discussed in Remark~\ref{rem:DiscreteOrthogonality}.

Restriction of our attention to the Fourier basis is motivated by its special property: discrete orthogonality. The next lemma illustrates this property for the sine basis (Example~\ref{exa:heat_equation}).

\begin{lem}[discrete orthogonality]
\label{lem:DiscreteOrthogonality}
Let $\{\psi_k\}_{k\in\mathbb{N}}$ be the sine basis, i.e.
\[
\psi_k (x) = \sqrt{2} \sin(k \pi x), \quad k = 1, 2, 3, \cdots.
\]
Then:
\begin{enumerate}[(i.)]
\item {Discrete orthogonality} holds:
\begin{equation}
\label{eq:DiscreteOrthogonality}
\langle \psi_j,\psi_k\rangle_d :=  \frac{1}{n}\sum_{i=1}^{n}\psi_j(i/n)\psi_k(i/n) = \delta_{jk},
\quad
j,k = 1,\cdots, n-1.
\end{equation}
Here $\delta_{jk}$ is the Kronecker delta.
\item Fix $l \in \mathbb{N}$. For any fixed $1 \leq k \le n-1$ and all $j\in\{ln,ln+1,\cdots,(l+1)n -1 \}$, there exits only one
$\bar{k}\in\{1,2,\ldots,n-1\}$ depending only on the parity of $l$, such that for $\tilde{j} = ln + \bar{k},$ the equality
\begin{equation}
\label{eq:DiscreteOrthogonalityTail}
| \langle \psi_{\tilde{j}},\psi_k\rangle_d | = 1
\end{equation}
holds, while $\langle \psi_{\tilde{j}},\psi_k\rangle_d =0 $ for all $\tilde{j} = ln + \tilde{k}$ such that $\tilde{k}\neq \bar{k},$ $\tilde{k}\in\{1,2,\ldots,n - 1\}.$
\end{enumerate}
\end{lem}

\begin{rem}\label{rem:DiscreteOrthogonality}
For other trigonometric bases, discrete orthogonality can also be attained.
Thus, the conjugate eigenbasis in Example \ref{exa:volterra_operator} is discretely orthogonal with design points $\{(i-1/2)/n\}_{i=1,\cdots,n}$. We refer to \cite{akansu2010GDFTs} and references therein for details. With some changes in the arguments, our asymptotic statistical results still remain valid with such modifications of design points compared to \eqref{eq:design_points}. We would like to stress the fact that restricting attention to bases with discrete orthogonality property does constitute a loss of generality. However, there exist classical bases other than trigonometric bases that are discretely orthogonal (possibly after a suitable modification of design points). See, for instance, \cite{quarteroni2010numerical} for an example of Lagrange polynomials.
\end{rem}

Motivated by the observations above, we introduce our central assumption on the basis functions.

\begin{ass}
\label{ass:discrete_orthogonality}
Given the design points $\{x_i\}_{i=1,\cdots, n}$ in \eqref{eq:design_points}, we assume the conjugate basis $\{\psi_k \}_{k\in\mathbb{N}}$ of the operator $A$ in \eqref{eq:general_inverse_problem_discrete} possesses the following properties:
\begin{enumerate}[(i.)]
\item for $1\leq j,k\leq n-1$,
\[
\langle \psi_j(x),\psi_k(x)\rangle_d:= \frac{1}{n} \sum_{i=1}^{n}\psi_j(x_i)\psi_k(x_i). 
=\delta_{jk}
\]

\item For $1\leq k \leq n-1$ and $j\in\{ln,\cdots,(l+1)n - 1\}$ with fixed $l \in \mathbb{N}$, there exits only one $\tilde{j} = ln + \bar{k},$ such that 
$
0 < | \langle \psi_{\tilde{j}},\psi_k\rangle_d | < M,
$
where $M$ is a fixed constant, and $\bar{k}$ depends on the parity of $l$ only.
For other $j \neq \tilde{j}$, $| \langle \psi_j,\psi_k\rangle_d | = 0.$
\end{enumerate}
\end{ass}

Using the shorthand notation
\[
f = \sum_j f_j \varphi_j = \sum_{j=1}^{n-1} f_j \varphi_j + \sum_{j \geq n} f_j \varphi_j =: f^n + f^r,
\]
we obtain for $k = 1,\cdots, n-1$ that
\begin{equation}
\label{eq:DisProj_Full}
\begin{split}
U_k =&  \frac{1}{n}\sum_{i=1}^n Y_i \psi_k(x_i) 
= \langle A f^n,\psi_k \rangle_d
	+\langle A f^r,\psi_k \rangle_d
	+ \frac{1}{n}\sum_{i=1}^n \xi_i \psi_k(x_i) \\
	=& a_k  f_k
	+ R_k
	+ \frac{1}{\sqrt{n}}\zeta_k,
\end{split}
\end{equation}
where
\[
R_k:= R_k(f) = \langle A f^r,\psi_k \rangle_d,\quad \zeta_k := \frac{1}{\sqrt{n}}\sum_{i=1}^n \xi_i \psi_k(x_i) .
\]
By Assumption \ref{ass:discrete_orthogonality}, we have 
\begin{equation}
\label{eq:RemainderTerm}
|R_k| = |\langle A f^r,\psi_k \rangle_d|
\leq\sum_{j \geq n} a_j |f_j| |\langle \psi_j,\psi_k\rangle_d|
= \sum_{l=1}^\infty a_{ln+\bar{k}} |f_{ln+\bar{k}}|,
\end{equation}
which leads to (via Cauchy-Schwarz)
\[
R_k^2(f)
\leq (\sum_{l=1}^\infty a_{ln+\bar{k}}^2 (ln+\bar{k})^{-2\beta}) \|f\|_{\beta}^2.
\]

Hence, for a mildly ill-posed problem, i.e. $a_k \asymp k^{-p},$ the following bound holds, uniformly in the ellipsoid $\{f: \|f\|_\beta \leq K \},$
\begin{align}
\label{eq:remainder_estimation_mild}
\sup_{f: \|f\|_\beta \le K} R_k^2(f) \lesssim & \sum_{l=1}^\infty (ln )^{-2\beta-2p} 
= n ^{-2(\beta+p)} \sum_{l=1}^\infty l^{-2(\beta+p)} \\
\asymp&  n ^{-2(\beta+p)} = o(1/n), \nonumber
\end{align}
for any $1 \leq k \leq n-1$ when $\beta+ p > 1/2$. 

If the problem is extremely ill-posed, i.e. $a_k \asymp e^{-k^sp},$ we use the inequality
\[
R_k^2(f)
\leq \left(\sum_{j \geq n} a_j|f_j| \right)^2 
\leq (\sum_{j \geq n} a_j^2) \|f^r\|^2.
\]
Since
$a_j \asymp \exp(-pj^s)\leq \exp(-pj)$, it follows that $\sum_{j \geq n}a_j^2$ is up to a constant bounded from above by $\exp(-2pn)$. Hence 
\begin{equation}
\label{eq:remainder_estimation_extreme}
\sup_{f: \|f\|_\beta \le K} R_k^2(f) \lesssim \exp(-2pn) \ll o(1/n).
\end{equation}

In \cite{knapik2011bayesianmild}, \cite{knapik2013bayesianextreme}, the Gaussian prior 
$\Pi = \otimes_{i\in\mathbb{N}} \mathcal{N}(0,\lambda_i)$ is employed on the coordinates of the eigenbasis expansion of $f$. If $\lambda_i = \rho_n^2 i^{-1-2\alpha},$ the sum $\sum_{i\in \mathbb{N}} \lambda_i = \rho_n^2 \sum_{i\in \mathbb{N}}i^{-1-2\alpha} $ is convergent, and hence this prior is the law of a Gaussian element in $H_1$.

In our case, we consider the same type of the prior with an additional constraint that only the first $n-1$ components of the prior are non-degenerate, i.e. 
$\Pi = \left(\otimes_{i < n}\mathcal{N}(0,\lambda_i)\right)\times\left(\otimes_{i\geq n}\mathcal{N}(0,0)\right)$, 
where $\lambda_i$ is as above. In addition, we assume the prior on $f$ is independent of the noise $\zeta_k, $ $k = 1,\cdots, n-1,$ in \eqref{eq:DisProj_Full}.
With these assumptions in force, we see $\Pi(R_k = 0) = 1,$ for $k = 1,\cdots, n-1.$ Furthermore, the posterior can be obtained from the product structure of the model and the prior via the normal conjugacy,
\begin{align}
\label{eq:gaussianInverse_posterior}
&\Pi(f | U^n) = \otimes_{k\in\mathbb{N}} \mathcal{N}(\hat{f}_k, \sigma_k^2), \\
&\text{with}\quad
\hat{f}_k = \frac{ n a_k \lambda_k \mathds{1}_{\{k < n\} } }{n a_k^2 \lambda_k + 1 } U_k, \quad
\sigma_k^2 = \frac{\lambda_k \mathds{1}_{\{k < n\} } } {n a_k^2 \lambda_k + 1}. \nonumber
\end{align}
We also introduce
\begin{align}
\label{eq:gaussianInverse_posteriorMean}
\hat{f} = \mathbb{E}(f | U^n ) = (
\mathbb{E}(f_k | U_k )) = ( \hat{f}_k )_{k\in\mathbb{N}} = (b_k U_k)_{k\in\mathbb{N}},
\end{align}
where $b_k = \frac{ n a_k \lambda_k \mathds{1}_{\{k < n\} } }{n a_k^2 \lambda_k + 1}.$
We conclude this section with a useful fact that will be applied in later sections:
\begin{align}
\label{eq:gaussianInverse_posteriorMean_coordinatewise_distr}
\hat{f}_k = b_k U_k = b_k \left( a_k f_k + R_k + \frac{\zeta_k}{\sqrt{n}}\right) = \mathbb{E} \hat{f}_k + \tau_k \zeta_k,
\end{align}
where $\mathbb{E} \hat{f}_k = a_k b_k f_k + b_k R_k $  and  $\tau_k = b_k/\sqrt{n}.$

\section{Main results}
\label{sec:MainResults}

\subsection{Contraction rates}

In this section, we determine the rate at which the posterior distribution concentrates on shrinking neighbourhoods of the `true' parameter $f_0$ as the sample size $n$ grows to infinity.

Assume the observations in \eqref{eq:general_inverse_problem_discrete} have been collected under the parameter value $f_0 = \sum_{k\in\mathbb{N}} f_{0,k} \varphi_k$. Thus our observations $(U_k)_{k<n}$ given in \eqref{eq:DisProj_Full} have the law $\otimes_{k < n}\mathcal{N} (a_k f_{0,k} + R_k, 1/n) $. We will use the notation $\Pi_n(\cdot | U)$ to denote the posterior distribution given in \eqref{eq:gaussianInverse_posterior}.

\begin{thm}[Posterior contraction: mildly ill-posed problem]
\label{thm:MildIllposed}
If the problem is mildly ill-posed as $a_k\asymp k^{-p}$ with $p > 0$, the true parameter $f_0 \in S^\beta$ with $\beta > 0$, and furthermore $\beta + p >1/2,$ by letting $\lambda_k = \rho_n^2 k ^{-1-2\alpha}$ with $\alpha>0$ and any positive $\rho_n$ satisfying $\rho_n^2 n\to \infty$,
we have, for any $K>0$ and $M_n\to \infty$,
\[
\sup_{\|f_0\|_\beta \leq K}
\mathbb{E}_{f_0} \Pi_n
\left(
f:\|f-f_0\|_{H_1}\geq M_n \varepsilon_n |U^n
\right)
\to 0,
\]
where
\begin{equation}
\label{eq:contraction_rate_mild}
\varepsilon_n 
=
\varepsilon_{n,1} \vee \varepsilon_{n,2}
=
(\rho_n^2 n)^{-\beta/(2\alpha+2p+1)\wedge 1}
\vee
\rho_n(\rho_n^2 n)^{-\alpha/(2\alpha+2p+1)}.
\end{equation}
In particular,
\begin{enumerate}[(i.)]
\item if $\rho_n = 1$, then $\varepsilon_n = n^{-(\alpha \wedge \beta)/(2\alpha+2p+1)}$;
\item if $ \beta \leq 2\alpha+2p+1$ and $\rho_n \asymp n^{(\alpha-\beta)/(2\beta +2p+1)}$, then $\varepsilon _n = n^{-\beta/(2\beta +2p+1)}$;
\item if $\beta > 2\alpha+2p+1$, then for every scaling $\rho_n$, $\varepsilon_n\gg n^{-\beta/(2\beta +2p+1)}.$
\end{enumerate}
\end{thm}

Thus we recover the same posterior contraction rates as obtained in \cite{knapik2011bayesianmild}, at the cost of an extra constraint $\beta + p > 1/2$. 
The frequentist minimax convergence rate for mildly ill-posed problems in the white noise setting with $\varepsilon=n^{-1/2}$ is $n^{-\beta/(2\beta + 2 p + 1)},$ see \cite{cavalier2008nonparametricinverseproblems}. We will compare our result to this rate. Our theorem states that in case (i.) the posterior contraction rate reaches the frequentist optimal rate if the regularity of the prior matches the truth $(\beta = \alpha)$ and the scaling factor $\rho_n$ is fixed. Alternatively, as in case (ii.), the optimal rate can also be attained by proper scaling, provided a sufficiently regular prior is used. In all other cases the contraction rate is slower than the minimax rate. Our results are similar to those in \cite{knapik2011bayesianmild} in the white noise setting. The extra constraint $\beta+p>1/2$ that we have in comparison to that work demands an explanation. 
As  \eqref{eq:remainder_estimation_mild} shows, the size of negligible terms $R_k(f_0)$ in \eqref{eq:DisProj_Full} decreases as the smoothness $\beta+p $ of the transformed signal $Af_0$ increases. In order to control $R_k$, a minimal smoothness of $Af_0$ is required. The latter is guaranteed if $p+\beta\geq 1/2,$ for it is known that in that case $Af_0$ will be at least continuous, while it may fail to be so if $p+\beta<1/2$, see \cite{tsybakov2008introduction}.

\begin{rem}
The control on $R_k(f_0)$ from \eqref{eq:RemainderTerm} depends on the fact that the eigenbasis possesses the properties in Assumption~\ref{ass:discrete_orthogonality}. If instead of Assumption~\ref{ass:discrete_orthogonality}~(ii.) one only assumes $|\langle \psi_j,\psi_k\rangle|\leq 1$ for any $k\leq n -1 $ and $j \geq n$, the constraint on the smoothness of $Af_0$ has to be strengthened to $\beta + p \geq 1$ in order to obtain the same results as in Theorem \ref{thm:MildIllposed}, because the condition $\beta + p \geq 1$ guarantees that the control on $R_k(f_0)$ in \eqref{eq:remainder_estimation_mild} remains valid.
\end{rem}

Now we consider the extremely ill-posed problem. The following result holds.

\begin{thm}[Posterior contraction: extremely ill-posed problem]
\label{thm:ExtremeIllposed}
Let the problem be extremely ill-posed as $a_k\asymp e^{-p k^s}$ with $s\geq 1$, and let the true parameter $f_0\in S^\beta$ with $\beta >0.$ Let $\lambda_k = \rho_n^2 k ^{-1-2\alpha}$ with $\alpha>0$ and any positive $\rho_n$ satisfying $\rho_n^2 n\to \infty$.
Then
\[
\sup_{\|f_0\|_\beta \leq K}
\mathbb{E}_{f_0} \Pi_n
\left(
f:\|f-f_0\|_{H_1}\geq M_n \varepsilon_n |U^n
\right)
\to 0,
\]
for any $K>0$ and $M_n\to \infty$, where
\begin{equation}
\label{eq:contraction_rate_extreme}
\varepsilon_n
=
\varepsilon_{n,1}\vee \varepsilon_{n,2}
= \left(\log(\rho_n^2 n)\right)^{-\beta/s} 
\vee
\rho_n\left(\log(\rho_n^2 n)\right)^{-\alpha/s}.
\end{equation}
In particular,
\begin{enumerate}[(i.)]
\item
if $\rho_n = 1$, then $\varepsilon_n = (\log n)^{-(\alpha\wedge\beta)/s}$,
\item
if $n^{-1/2+\delta}\lesssim \rho_n \lesssim(\log n)^{(\alpha - \beta)/s}$ for some $\delta >0$, then $\varepsilon_n = (\log n)^{-\beta/s}$.
\end{enumerate}

Furthermore, if $\lambda_k = \exp(-\alpha k^s)$ with $\alpha >0$, the following contraction rate is obtained:
$\varepsilon_n=(\log n)^{-\beta/s}.$
\end{thm}
Since the frequentist minimax estimation rate in extremely ill-posed problems in the white noise setting is $(\log n)^{-\beta/s}$ (see \cite{cavalier2008nonparametricinverseproblems}), Theorem~\ref{thm:ExtremeIllposed} shows that the optimal contraction rates can be reached by suitable choice of the regularity of the prior, or by using an appropriate scaling. In contrast to the mildly ill-posed case, we have no extra requirement on the smoothness of $Af_0$. The reason is obvious: because the signal is lifted to $S^\infty$ by the forward operator $A$, the term \eqref{eq:remainder_estimation_extreme} converges to zero exponentially fast, implying that $R_k(f_0)$ in \eqref{eq:DisProj_Full} is always negligible.

\subsection{Credible sets}
\label{sec:credible_sets}

In the Bayesian paradigm, the spread of the posterior distribution is a common measure of uncertainty in parameter estimates. In this section we study the frequentist coverage of Bayesian credible sets in our problem.

When the posterior is Gaussian, it is customary to consider credible sets centered at the posterior mean, which is what we will also do. In addition, because in our case the covariance operator of the posterior distribution does not depend on the data, the radius of the credible ball is determined by the credibility level $1-\gamma$ and the sample size $n$. A credible ball centred at the posterior mean $\hat{f}$ from \eqref{eq:gaussianInverse_posteriorMean} is given by
\begin{equation}
\label{eq:credible_region}
\hat{f} + B(r_{n,\gamma}) := \{f\in H_1: \|f- \hat{f} \|_{H_1} \leq r_{n,\gamma} \},
\end{equation}
where the radius $r_{n,\gamma}$ is determined by the requirement that
\begin{equation}
\label{eq:r_n}
\Pi_n (\hat{f} + B(r_{n,\gamma}) | U^n) = 1 - \gamma.
\end{equation}

By definition, the frequentist {coverage} or confidence of the set \eqref{eq:credible_region} is
\begin{equation}
\label{eq:frequentist_coverage}
\mathbb{P}_{f_0} (f_0 \in \hat{f} + B(r_{n,\gamma}) ),
\end{equation}
where the probability measure is the one induced by the law of $U^n$ given in \eqref{eq:DisProj_Full} with $f = f_0$. We are interested in the asymptotic behaviour of the coverage \eqref{eq:frequentist_coverage} as $n\to \infty$ for a fixed $f_0$ uniformly in Sobolev balls, and also along a sequence $f_0^n$ changing with $n$.

The following two theorems hold.

\begin{thm}[Credible sets: mildly ill-posed problem]
\label{thm:credibleset_mildly}
Assume the same assumptions as in Theorem~\ref{thm:MildIllposed} hold, and let $\tilde{\beta} = \beta \wedge (2\alpha + 2 p +1)$. The asymptotic coverage of the credible set \eqref{eq:credible_region} is
\begin{enumerate}[(i.)]
\item
1, uniformly in $\{f_0: \|f_0\|_{\beta}\leq 1\}$, if $\rho_n \gg n^{(\alpha- \tilde{\beta})/(2\tilde{\beta}+2p+1)}$;

\item
1, for every fixed $f_0\in S^\beta$, if $\beta < 2\alpha + 2p +1$ and $\rho_n \asymp n^{(\alpha-\tilde{\beta})/(2\tilde{\beta}+2p+1)}$;
c, along some $f_0^n$ with $\sup_n \|f_0^n\|_{\beta} <\infty$, if $\rho_n \asymp n^{(\alpha-\tilde{\beta})/(2\tilde{\beta}+2p+1)}$ (any $c\in[0,1)$).

\item
0, along some $f_0^n$ with $\sup_n \|f_0^n\|_{\beta} <\infty$, if $\rho_n \ll n^{(\alpha-\tilde{\beta})/(2\tilde{\beta}+2p+1)}$.
\end{enumerate}
\end{thm}

\begin{thm}[Credible sets: extremely ill-posed problem]
\label{thm:credibleset_extremely}
Assume the setup of Theorem~\ref{thm:ExtremeIllposed}. Then if $\lambda_k = \rho_n^2 k ^{-1-2\alpha}$ with $\alpha>0$ and any positive $\rho_n$ satisfying $\rho_n^2 n\to \infty$, the asymptotic coverage of the credible set \eqref{eq:credible_region} is
\begin{enumerate}[(i.)]
\item
1, uniformly in $\{f_0: \|f_0\|_{S^\beta}\leq 1\}$, if $\rho_n \gg (\log n)^{(\alpha-\beta)/2}$;

\item
1, uniformly in $f_0$ with $\|f_0\|_{\beta}\leq r$ with $r$ small enough;

1, for any fixed $f_0\in S^\beta$,\\
provided the condition $\rho_n \asymp (\log n)^{(\alpha-\beta)/s}$ holds;
\item
0, along some $f_0^n$ with $\sup_n \|f_0^n\|_{\beta} <\infty$, if $\rho_n \lesssim (\log n)^{(\alpha - \beta)/s}$.

\newcounter{enumTemp}
\setcounter{enumTemp}{\value{enumi}}
\end{enumerate}

Moreover, if $\lambda_k = e^{-\alpha^s}$ with $\alpha>0$ and any positive $\rho_n$ satisfying $\rho_n^2 n\to \infty$, the asymptotic coverage of the credible set \eqref{eq:credible_region} is
\begin{enumerate}[(i.)]
    \setcounter{enumi}{\theenumTemp}
\item
0, for every $f_0$ such that $|f_{0,i}|\gtrsim e^{-ci^s/2}$ for some $c <\alpha$.
\end{enumerate}
\end{thm}

For the two theorems in this section, the most intuitive explanation is offered by the case $\rho_n \equiv 1$. The situations (i.), (ii.) and (iii.) correspond to $\alpha<\beta$, $\alpha= \beta$ and $\alpha>\beta$, respectively. The message is that the oversmoothing prior ((iii.) in Theorem~\ref{thm:credibleset_mildly} and (iii.), (iv.) in Theorem~\ref{thm:credibleset_extremely}) leads to disastrous frequentist coverage of credible sets, while the undersmoothing prior ((i.) in both theorems) delivers very conservative frequentist results (coverage 1). With the right regularity of the prior (case (ii.)), the outcome depends on the norm of the true parameter $f_0$. Our results are thus similar to those obtained in the white noise setting in \cite{knapik2011bayesianmild} and \cite{knapik2013bayesianextreme}.

\section{Simulation examples}
\label{sec:SimulationExamples}

In this section we carry out a small-scale simulation study illustrating our theoretical results. Examples we use to that end are those given in Subsection~\ref{subsec:SVD}. These were also used in simulations in \cite{knapik2011bayesianmild} and \cite{knapik2013bayesianextreme}.

In the setting of Example~\ref{exa:volterra_operator}, we use the following true signal,
\begin{equation}
\label{eq:volterra_euqation_true_parameter}
f_0(x) = \sum_{i=1}^{\infty} f_{0,i} \varphi_i(x) \text{ with } f_{0,k} = k^{-3/2} \sin (k).
\end{equation}
It is easy to check that $f_0 \in S^1$.

In the setup of Example~\ref{exa:heat_equation}, the initial condition is assumed to be
\begin{equation}
\label{eq:heat_equation_initial_condition}
f_0(x) = 4x(x-1)(8x-5).
\end{equation}
One can verify that in this case
\begin{equation*}
f_{0,k} = \frac{8\sqrt{2}(13+11(-1)^k)}{\pi^3 k^3},
\end{equation*}
and $f_0 \in S^\beta$ for any $\beta < 5/2$.

First, we generate noisy observations $\{Y_i\}_{i=1,\cdots,n}$ from our observation scheme \eqref{eq:general_inverse_problem_discrete} at design points $x_i = \frac{i-1/2}{n}$ in the case of Volterra operator, and $x_i = i/n$ in the case of the heat equation. Next, we apply the transform described in \eqref{eq:DisProj_Full} and obtain transformed observations $\{U_i\}_{i=1,\cdots,n-1}$. Then, by \eqref{eq:gaussianInverse_posterior}, the posterior of the coefficients with the eigenbasis $\varphi_i$ is given by
\[
f_k | U^n \sim \mathcal{N}\left(\frac{n a_k \lambda_k \mathds{1}_{\{k<n\} } }{n a_k^2\lambda_k +1} U_k,
				\frac{\lambda_k \mathds{1}_{\{k<n\} } }{n a_k^2\lambda_k +1} \right).
\]

\begin{figure}
\centering
\includegraphics[width=0.8\textwidth]{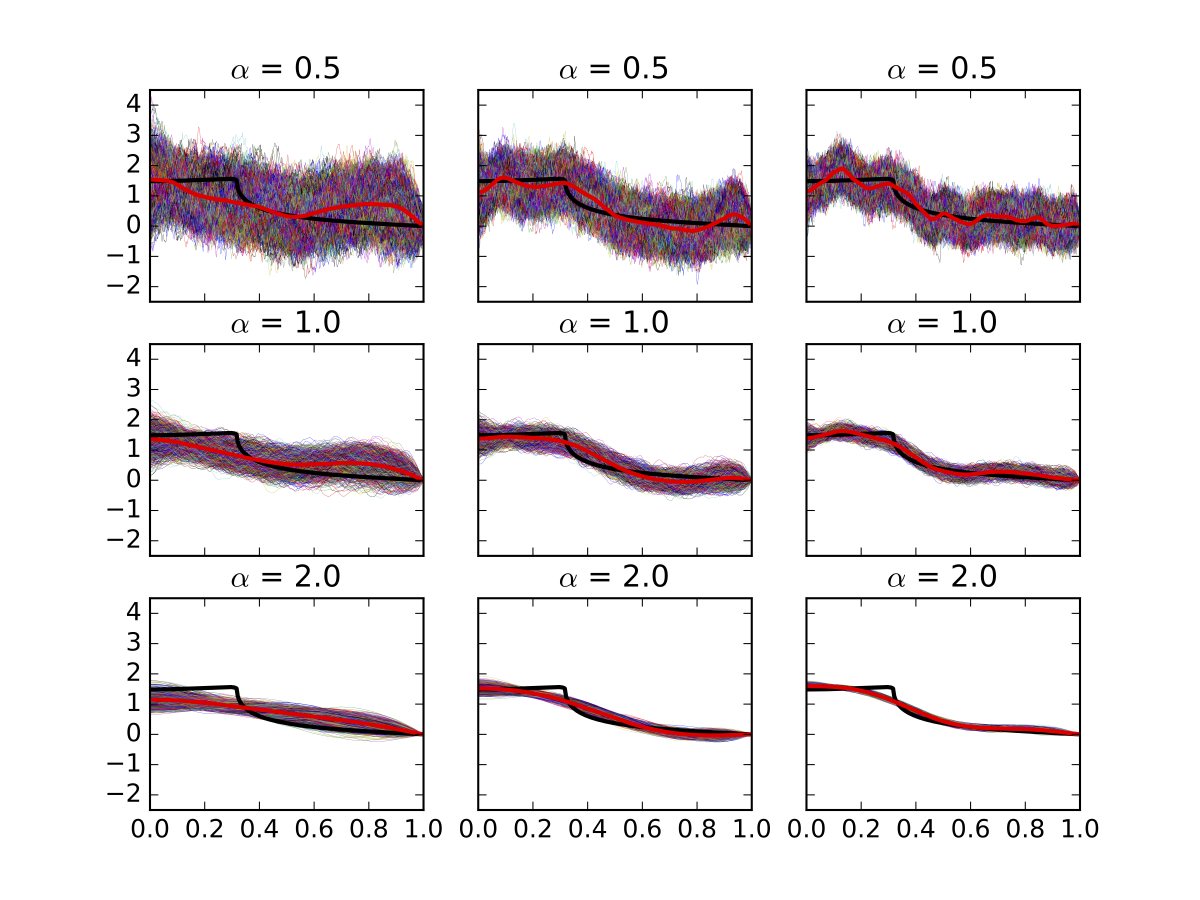}
\caption{Realizations of the posterior mean (red) and 950 of 1000 draws from the posterior (colored thin lines) with smallest $L^2$ distance to the posterior mean. From left to right columns, the posterior is computed based on sample size $10^3, 10^4$ and $10^5$ respectively. The true parameter (black) is of smoothness $\beta = 1$ and given by coefficients $f_{0,k} = k^{-3/2} \sin(k)$.}
\label{fig:volterra_operator}
\end{figure}

\begin{figure}
\centering
\includegraphics[width=0.8\textwidth]{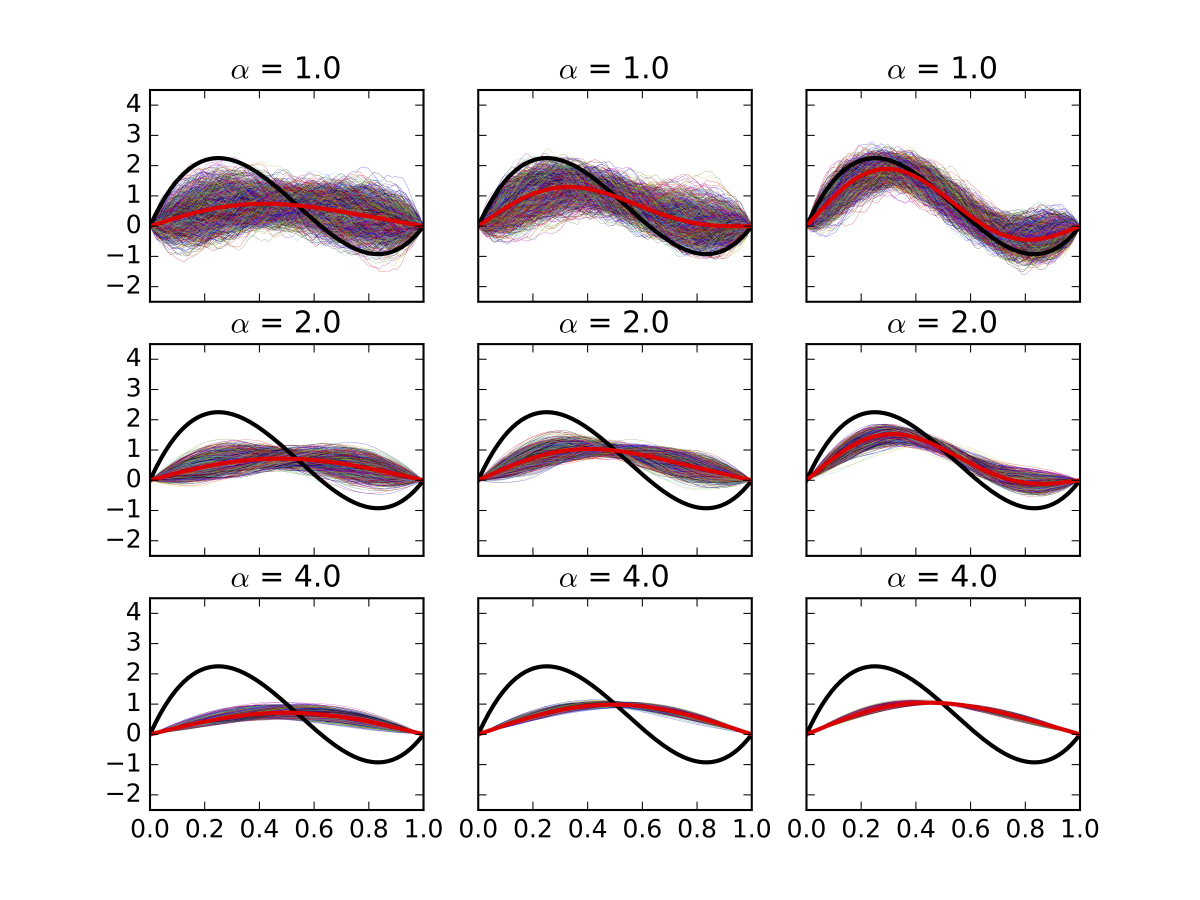}
\caption{Realizations of the posterior mean (red) and 950 of 1000 draws from the posterior (colored thin lines) with smallest $L^2$ distance to the posterior mean. From left to right columns, the posterior is computed based on sample size $10^3, 10^4$ and $10^5$ respectively. The true parameter (black) is of smoothness $\beta$ for any $\beta<5/2$ and given by  \eqref{eq:heat_equation_initial_condition}.}
\label{fig:heat_equation}
\end{figure}

Figures~\ref{fig:volterra_operator} and \ref{fig:heat_equation} display plots of $95\%$ $L_2$-credible bands for different sample sizes and different priors. For all priors we assume $\rho_n \equiv 1$, and use different smoothness degrees $\alpha$, as shown in the titles of the subplots. In addition, the columns from left to right corresponds to $10^3, 10^4$ and $10^5$ observations. The (estimated) credible bands are obtained by generating $1000$ realizations from the posterior and retaining $95\%$ of them that are closest in the $L^2$-distance to the posterior mean.


Two simulations reflect several similar facts. First, because of the difficulty due to the inverse nature of the problem, the recovery of the true signal is relatively slow, as the posteriors for the sample size $10^3$ are still rather diffuse around the true parameter value. Second, it is evident that undersmoothing priors (the top rows in the figures) deliver conservative credible bands, but still capture the truth. On the other hand, oversmoothing priors lead to overconfident, narrow bands, failing to actually express the truth (bottom rows in the figures). As already anticipated due to a greater degree of ill-posedness, recovery of the initial condition in the heat equation case is more difficult than recovery of the true function in the case of the Volterra operator. Finally, we remark that qualitative behaviour of the posterior in our examples is similar to the one observed in \cite{knapik2011bayesianmild} and \cite{knapik2013bayesianextreme}; for larger samples sizes $n$, discreteness of the observation scheme does not appear to have a noticeably adversary effect compared to the fully observed case in \cite{knapik2011bayesianmild} and \cite{knapik2013bayesianextreme}.

\section{Proofs}
\label{sec:proofs}

\subsection{Proof of Lemma~\ref{lem:DiscreteOrthogonality}}

This proof is a modification of the one of Lemma 1.7 in \cite{tsybakov2008introduction}. With the following temporary definitions $a:= e^{i \pi \frac{j}{n}}$ and $\quad b:= e^{i \pi \frac{k}{n}},$ using Euler's formula, we have
\begin{equation}
\label{eq:product}
\begin{split}
\langle \psi_j,\psi_k\rangle_d =& 
-\frac{1}{2n}
\sum_{s=1}^{n}
(a^s - a^{-s})(b^s - b^{-s})\\
=&
-\frac{1}{2n}
\sum_{s=1}^{n}
\left[
(ab)^s - (a/b)^s - (a/b)^{-s} + (ab)^{-s}
\right],\\
=&
-\frac{1}{2n}
\left[
\underbrace{\sum_{s=1}^{n}(ab)^s}_A - \underbrace{\sum_{s=1}^{n}(a/b)^s}_B - \underbrace{\sum_{s=1}^{n}(a/b)^{-s}}_C + \underbrace{\sum_{s=1}^{n}(ab)^{-s}}_D
\right].
\end{split}
\end{equation}
Furthermore,
\[
ab = e^{i\pi \frac{j+k}{n}}, \ \frac{a}{b} = e^{i\pi \frac{j-k}{n}}.
\]
Observe that when $ab \neq 1$, we have 
\begin{align*}
A = \frac{ab(1-(ab)^n)}{1 - ab}, \quad D = \frac{1 - (ab)^{-n}}{ab - 1}, \quad
A + D = \frac{ab(1-(ab)^n) - (1 - (ab)^{-n})}{1 - ab}.
\end{align*}
Similarly, if $a/b \neq 1$,
\begin{align*}
B + C = \frac{(a/b)(1-(a/b)^n) - (1 - (a/b)^{-n})}{1 - (a/b)}.
\end{align*}

We fix $1 \leq k \leq n-1$ and discuss different situations depending on $j$. 

\begin{enumerate}[(I.)]
\item $1 \leq j \leq n-1$ and $j + k \neq n$.

Since $n \neq j+k < 2n$, we always have $ab = e^{i\pi \frac{j+k}{n}}\neq 1$, and the terms $A$ and $D$ can be calculated as above. Similarly, since $-n < j - k < n$, $a/b = 1$ only when $j = k$. Moreover, $j+k$ and $j-k$ have the same parity, and so $j=k$ is only possible if $j+k$ is even.
\begin{enumerate}[(i.)]
\item $j + k$ is even. 

In this case, $(ab)^n = 1$. This leads to $A = D = 0$.

Further, if $j = k$, we have $a/b = b/a = 1$ and $B = C = n$. Otherwise, if $j\neq k$, we have $a/b \neq 1$ and $(a/b)^n = 1 = (b/a)^n$ (since $j-k$ is even), and so
\begin{align*}
B = \frac{a/b(1-(a/b)^n)}{1-a/b} = 0, \quad C = 0,
\end{align*}
which implies \eqref{eq:product} equals $1$.

\item $j + k$ is odd.
We have $(ab)^{n}  = (a/b)^{n} = -1$, which results in $A + D = B + C = -2$, and so \eqref{eq:product} equals $0$.

\end{enumerate}

\item $1 \leq j <n$ and $j + k = n$.
We have $ab = -1$. Arguing as above, if $n$ is odd, $A + D = -2$ and $B + C = -2$. If $n$ is even, $A = D = 0$ and $B = C = n \delta_{jk}$.
\setcounter{enum}{\value{enumi}}
\end{enumerate}

The remaining cases follow the same arguments, and hence we omit the (lengthy and elementary) calculations.

\begin{enumerate}[(I.)]
\setcounter{enumi}{\value{enum}}
\item $j = ln$ with $l \in\mathbb{N}$. 

It can be shown that $A+D = B+C$ always holds.

\item $j \in\{ln+1,\cdots,(l+1)n-1\}$. 

When $l$ is even, one obtains $\langle \psi_j,\psi_k\rangle_d = \delta_{\tilde{j}k}$, where $\tilde{j} = j - ln$. Otherwise, for odd $l$, $\langle \psi_j,\psi_k\rangle_d = -\delta_{\tilde{j}k}$ where $\tilde{j} = (l+1)n - j$.
\end{enumerate}


\subsection{Proof of Theorem~\ref{thm:MildIllposed}}
\label{sec:proof_mild}


In this proof we use the notation $\|\cdot\| = \|\cdot\|_{H_1} = \|\cdot\|_{\ell^2}$. To show
\begin{align*}
\sup_{\|f_0\|_\beta\leq K}\mathbb{E}_{f_0}\Pi_n\left(
f:\|f-f_0\|\geq M_n \varepsilon_n |U^n
\right)
\to 0,
\end{align*}
we first apply Markov's inequality,
\begin{align*}
 M_n^2\varepsilon_{n}^2 \Pi_n \left( f: \|f-f_0\|^2\geq M_n^2 \varepsilon_{n}^2 | U^n\right)
\leq \int \|f-f_0\|^2 \,\mathrm{d}\Pi_n(f|U^n).
\end{align*}
From \eqref{eq:gaussianInverse_posterior} and the bias-variance decomposition,
\begin{align*}
\int \|f-f_0\|^2 \,\mathrm{d}\Pi_n(f|U^n)
=\|\hat{f}-f_0\|^2 + \|\sigma\|^2,
\end{align*}
where $\sigma = (\sigma_k)_k$ is given in \eqref{eq:gaussianInverse_posterior}. Because $\sigma$ is deterministic,
\begin{align*}
\mathbb{E}_{f_0}\left[ 
	\Pi_n \left( f : \|f-f_0\|\geq M_n \varepsilon_n | U^n\right)
\right]
\leq
\frac{1}{M_n^2\varepsilon_n^2}
\left(
\mathbb{E}_{f_0} \|\hat{f} -f_0\|^2 + \|\sigma\|^2
\right).
\end{align*}
Since $M_n\to \infty$ is assumed, it suffices to show that the terms in brackets are bounded by a constant multiple of $\varepsilon_n^2$ uniformly in $f_0$ in the Sobolev ellipsoid. 


Using \eqref{eq:gaussianInverse_posteriorMean_coordinatewise_distr}, we obtain
\begin{align*}
\mathbb{E}_{f_0} \|\hat{f} -f_0\|^2
=
\|\mathbb{E}_{f_0}\hat{f}  -f_0\|^2 + \|\tau\|^2
=\|\mathbb{E}_{f_0}\hat{f}-f_0^n\|^2 + \|f_0^r\|^2 + \|\tau\|^2,
\end{align*}
where $\tau = (\tau_k)_k$ given in \eqref{eq:gaussianInverse_posteriorMean_coordinatewise_distr} and
\begin{align*}
 f_0^{n} =& (f_{0,1},\cdots,f_{0,n-1},0,\cdots),\\
  f_0^{r} =& (0,\cdots,0, f_{0,n}, f_{0,n+2},\cdots).
\end{align*}
We need to obtain a uniform upper bound over the ellipsoid $\{f_0: \|f_0\|_{\beta} \leq K\}$ for
\begin{align}
\label{eq:boundsdecomposition}
\|\mathbb{E}_{f_0}\hat{f}  -f_0^n\|^2 +\|f_0^r\|^2 + \|\tau\|^2 + \|\sigma\|^2.
\end{align}
We have
\begin{align}
\label{eq:square_bias}
&\|\mathbb{E}_{f_0}\hat{f}  -f_0^n\|^2 
=
\sum_{k = 1}^{n-1} 
\left(
\frac{n a_k^2\lambda_k}{n a_k^2\lambda_k+1} f_{0,k}
+
\frac{n a_k\lambda_k}{n a_k^2\lambda_k+1} R_k
-
f_{0,k}
\right)^2\nonumber\\
\lesssim&
\underbrace{
\sum_{k = 1}^{n-1} 
\frac{1}{\left(n a_k^2\lambda_k+1\right)^2}
f_{0,k}^2
}_{A_1}
+
n \sup_{k < n} R_k^2
\underbrace{
\sum_{k = 1}^{n-1} 
\frac{n  a_k^2\lambda_k^2}{\left(n a_k^2\lambda_k+1\right)^2}
}_{A_2},
\end{align}
and
\[
\|f_0^r\|^2 = \sum_{k \geq n} f_{0,k}^2,
\quad
\|\tau\|^2
=
\sum_{k=1}^{n-1}
\frac{n a_k^2\lambda_k^2}{\left(n a_k^2\lambda_k+1\right)^2} = A_2,  
\quad
\|\sigma\|^2
=
\sum_{k=1}^{n-1}
\frac{\lambda_k}{n a_k^2\lambda_k+1}.
\]
Recall that we write \eqref{eq:contraction_rate_mild} as $\varepsilon_n = \varepsilon_{n,1} \vee \varepsilon_{n,2}$. The statements (i.)--(iii.) follow by elementary calculations. Specifically, in (ii.) the given $\rho_n$ is the best scaling, as it gives the fastest rate. From \cite{knapik2011bayesianmild} (see the argument below (7.3) on page~21), $A_1$ is bounded by a fixed multiple of $(\varepsilon_{n,1})^2,$ and $\|\tau\|^2$, $\|\sigma\|^2$ are bounded by multiples of $(\varepsilon_{n,2})^2$. Hence, to show that the rate is indeed \eqref{eq:contraction_rate_mild}, it suffices to show that
$
n \sup_{k\leq n} R_k^2
A_2
$
and
$
\|f_0^r\|^2
$
can be bounded by a multiple of $(\varepsilon_n)^2$ uniformly in the ellipsoid $\{f_0: \|f_0\|_{\beta} \leq K\}$.
Since $A_2 = \|\tau\|^2$, to that end it is sufficient to show that $\sup_{k < n} n R_k^2 = O(1),$ and that $\|f_0^r\|^2=O(\varepsilon_n)^2.$

Since $f_0\in S^\beta$, we have the following straightforward bound,
\[
\|f_0^r\|^2 \leq n^{-2\beta}\sum_{k \geq n} f_{0,k}^2 k^{2\beta} \leq n^{-2\beta} \|f_0\|_{\beta}^2 \lesssim n^{-2\beta},
\]
which is uniform in $\{f_0: \|f_0\|_{\beta} \leq K\}.$ By comparing to the rates in the statements (ii.)--(iii.), it is easy to see that $n^{-2\beta}$ is always negligible with respect to $\varepsilon_n^2$.

Proving $\sup_{k\leq n}n R_k^2 = O(1)$ is equivalent to showing $\sup_{k\leq n}R_k^2 = O(1/n)$; but the latter has been already proved in \eqref{eq:remainder_estimation_mild}. Notice that we actually obtained a sharper bound $\sup_{k\leq n}n R_k^2=o(1)$ than the one necessary for our purposes in this proof. However, this sharper bound will be used in the proof of Theorem~\ref{thm:credibleset_mildly}. By taking supremum over $f_0$, we thus have
\begin{equation}
\label{eq:squrebias_mild_upperbound}
\sup_{\|f_0\|_{S^\beta}\leq K}  \left( \|\mathbb{E}_{f_0}\hat{f}  -f_0^n\|^2 +\|f_0^r\|^2 \right)
\lesssim \varepsilon_n^2 + n^{-2\beta}\lesssim \varepsilon_n^2,
\end{equation}
with which we conclude that up to a multiplicative constant, \eqref{eq:boundsdecomposition} is bounded by $\varepsilon_n^2$  uniformly over the ellipsoid $\sup_{\|f_0\|_{\beta}\leq K}$. This completes the proof.

\subsection{Proof of Theorem~\ref{thm:ExtremeIllposed}}
\label{sec:proof_extreme}
We start by generalizing Theorem 3.1 in \cite{knapik2013bayesianextreme}. Following the same lines as in the proof of that theorem and using Lemma~\ref{lem:supremum_estimation_extreme}, \ref{lem:A4}, \ref{lem:A5}, \ref{lem:A6} in Section~\ref{sec:appendix} of the present paper instead of analogous technical results in \cite{knapik2013bayesianextreme}, the statement of Theorem 3.1 in \cite{knapik2013bayesianextreme} can be extended from $s=2$ to a general $s\geq 1$, for which the posterior rate is given by \eqref{eq:contraction_rate_extreme}, or $\varepsilon_n = \varepsilon_{n,1} \vee \varepsilon_{n,2}$ in short.

In our model, we again obtain \eqref{eq:boundsdecomposition} and also that a fixed multiple of $(\varepsilon_{n,1})^2$ is an upper bound of $A_1,$ and that $\|\tau\|^2, \|\sigma\|^2$ can be bounded from above by fixed multiples of $(\varepsilon_{n,1})^2$.

Now as in the proof of Theorem~\ref{thm:MildIllposed} in Section~\ref{sec:proof_mild}, we will show that $\sup_{\|f_0\|_{\beta}\le K}(\|\mathbb{E}_{f_0}\hat{f}  -f_0^n\|^2 + \|f_0^r\|^2)$ can be bounded by a fixed multiple of $(\varepsilon_n)^2$ by proving that $\sup_{k\leq n} n R_k^2 = O(1)$. By \eqref{eq:remainder_estimation_extreme}, $n (R_k)^2 \lesssim \exp(-2pn)n,$ and the righthand side converges to zero. Therefore,
\[
\sup_{\|f_0\|_{\beta}\leq K} \left( \|\mathbb{E}_{f_0}\hat{f}  -f_0^n\|^2 +\|f_0^r\|^2 \right)
\lesssim \varepsilon_n.
\]
Parts (i.) and (ii.) of the statement of the theorem are obtained by direct substitutions, using the fact that $\log n \ll n$. Notice that if $\rho_n \gtrsim (\log n)^{(\alpha - \beta)/s}$, the rate $\varepsilon_n$ deteriorates and is dominated by the second term in \eqref{eq:contraction_rate_extreme}.

For the case $\lambda_k = \exp(-\alpha k^s)$, the argument follows the same lines as in Section 5.1 in \cite{knapik2013bayesianextreme}, and our arguments above.

\subsection{Proof of Theorem~\ref{thm:credibleset_mildly}}
\label{sec:proof_credibleset_mild}

The proof runs along the same lines as the proof of Theorem 4.2 in \cite{knapik2011bayesianmild}. We will only show the main steps here.

In Section \ref{subsec:gaussianInverse_equivalentFormulation}, we have shown that the posterior distribution is $\otimes_{k\in\mathbb{N}} \mathcal{N}(\hat{f}_k, \sigma_k^2)$, the radius $r_{n,\gamma}$ in \eqref{eq:credible_region} satisfies $\mathbb{P}_{X_n}(X_n <r_{n,\gamma}^2) = 1 - \gamma$, where $X_n$ is a random variable distributed as the square norm of an $\otimes_{k\in\mathbb{N}} \mathcal{N}(\hat{f}_k, \sigma_k^2)$ variable. Let $T = (\tau_k^2)_{k\in\mathbb{N}}.$ Under \eqref{eq:DisProj_Full}, the variable $\hat{f}$ is distributed as $\mathcal{N}_{H_1}(\mathbb{E}_{f_0}\hat{f}, T) \coloneqq  \otimes_{k\in\mathbb{N}}\mathcal{N}(\mathbb{E}_{f_0}\hat{f}_k, \tau_k^2).$  Hence the coverage \eqref{eq:frequentist_coverage} can be rewritten as 
\begin{equation}
\label{eq:frequentist_coverage_shifted}
\mathbb{P}_{W_n}(\|W_n +\mathbb{E}_{f_0}\hat{f} - f_0\|_{H_1} \leq r_{n,\gamma}),
\end{equation}
where $W_n \sim \mathcal{N}_{H_1}(0,T).$ Denote $V_n = \|W_n\|_{H_1}^2$ and observe that one has in distribution
\[
X_n = \sum_{1\leq i< n} \sigma^2_{i} Z_i^2,\quad V_n = \sum_{1\leq i < n} \tau^2_{i} Z_i^2
\]
for $\{Z_i\}$ independent standard Gaussian random variables  
with
\[
\sigma^2_{i} = \frac{\lambda_i}{n a_i^2\lambda_i + 1},
\quad \tau^2_{i} = \frac{n a_i^2 \lambda_i^2}{(n a_i^2\lambda_i + 1)^2}.
\]

By the same argument as in \cite{knapik2011bayesianmild}, one can show that the standard deviations of $X_n$ and $V_n$ are negligible with respect to their means,
\begin{equation}
	\label{eq:means}
\mathbb{E} X_n \asymp \rho_n^2(\rho_n^2 n)^{-2\alpha/(2\alpha + 2 p +1)},
\quad
\mathbb{E} V_n \asymp \rho_n^2(\rho_n^2 n)^{-2\alpha/(2\alpha + 2 p +1)},
\end{equation}
and the difference of their means,
\[
\mathbb{E}(X_n - V_n) \asymp \rho_n^2(\rho_n^2 n)^{-2\alpha/(2\alpha + 2 p +1)}.
\]
Since $X_n\geq V_n$, the distributions of $X_n$ and $V_n$ are asymptotically separated, i.e. $\mathbb{P}(V_n \leq v_n \leq X_n) \to 1$ for some $v_n$, e.g. $v_n = \mathbb{E}(V_n+X_n)/2$. Since $r_{n,\gamma}^2$ are $1-\gamma$ quantiles of $X_n$, we also have $\mathbb{P}(V_n\leq r_{n,\gamma}^2 (1+o(1))) \to 1$. In addition, by \eqref{eq:means},
\[
r_{n,\gamma}^2 \asymp \rho_n^2(\rho_n^2 n)^{-2\alpha/(2\alpha + 2 p +1)}.
\]
Introduce
\begin{equation}
\label{eq:square_norm_of_bias}
B_n:= \sup_{\|f_0\|_{\beta}\lesssim 1} \|\mathbb{E}_{f_0} \hat{f} - f_0\|_{H_1} = \sup_{\|f_0\|_{\beta}\lesssim 1} \left( \|\mathbb{E}_{f_0}\hat{f}  -f_0^n\|_{H_1} +\|f_0^r\|_{H_1} \right).
\end{equation}
It follows from the arguments for \eqref{eq:remainder_estimation_mild} in the proof of Theorem~\ref{thm:MildIllposed} that
\[
B_n \lesssim \varepsilon_{n,1} \vee \left(\sqrt{n} R \varepsilon_{n,2}\right),
\]
where $R = \sup_{k < n} R_k \lesssim n^{-(p+\beta)}$. Now apply the argument on the lower bound from Lemma 8.1 in \cite{knapik2011bayesianmild} (with $q =\beta, t=0, u = 2\alpha+2p+1, v = 2, N = \rho_n^2 n$) to obtain that $B_n \gtrsim \varepsilon_{n,1}$. Thus we have
\[
\varepsilon_{n,1}
\lesssim
B_n
\lesssim
\varepsilon_{n,1}\vee\left(\sqrt{n}R \varepsilon_{n,2}\right).
\]

We consider separate cases. In case (i.), substituting the corresponding $\rho_n$ into the expression of $\varepsilon_{n,1}$ and $\varepsilon_{n,2}$, we have $\varepsilon_{n,1}\ll \varepsilon_{n,2}$. By \eqref{eq:remainder_estimation_mild}, $B_n\lesssim\varepsilon_{n,1}\vee\left(\sqrt{n}R \varepsilon_{n,2}\right) \ll \varepsilon_{n,2} \asymp r_{n,\gamma}$. This leads to 
\begin{align}
\label{eq:coverage_B_n_ll_r_n}
\mathbb{P}(\|W_n +\mathbb{E}_{f_0} \hat{f} - f_0 \|_{H_1}\leq r_{n,\gamma})
\geq& 
\mathbb{P}(\|W_n \|_{H_1}\leq r_{n,\gamma} - B_n)\nonumber\\
=&
\mathbb{P}(V_n\leq r_{n,\gamma}^2(1+o(1))) \to 1
\end{align}
uniformly in the set $\{f_0: \|f_0\|_{\beta}\lesssim 1\}$.

In case (iii.), the given $\rho_n$ leads to $\varepsilon_{n,1}\gg \varepsilon_{n,2}$ and consequently $B_n\gg r_{n,\gamma}$. Hence,
\begin{align*}
\mathbb{P}(\|W_n +\mathbb{E}_{f_0} \hat{f}^n - f_0^n \|_{H_1}\leq r_{n,\gamma})
\leq
\mathbb{P}(\|W_n \|_{H_1}\geq B_n - r_{n,\gamma})\nonumber \to 0,
\end{align*}
for any $f_0^n$ (nearly) attaining the supremum.

In case (ii.), we have $B_n \asymp r_{n,\gamma}$. If $\beta < 2\alpha + 2 p +1$, by Lemma 8.1 in \cite{knapik2011bayesianmild} the bias $\mathbb{E}_{f_0} \hat{f} - f_0$ at a fixed $f_0$ is of strictly smaller order than $B_n$. Following the argument of case (i.), the asymptotic coverage can be shown to converge to 1.

For existence of a sequence along which the coverage is $c\in[0,1)$, we only give a sketch of the proof here; the details can be filled in as in \cite{knapik2011bayesianmild}.

The coverage \eqref{eq:frequentist_coverage_shifted} with $f_0$ replaced by $f_0^n$ tends to $c$, if for $b_n = \mathbb{E}_{f_0} \hat{f}^n - f_0^n$ and $z_c$ a standard normal quantile,
\begin{align}
\label{eq:existence_normal}
\frac{\|W_n+b_n\|_{H_1}^2 - \mathbb{E} \|W_n+b_n\|_{H_1}^2}{\operatorname{sd}\|W_n+b_n\|_{H_1}^2}
&\rightsquigarrow \mathcal{N}(0,1),\\
\label{eq:existence_quantile}
\frac{r_{n,\gamma}^2 - \mathbb{E}\|W_n + b_n\|_{H_1}^2}{\operatorname{sd}\|W_n + b_n\|_{H_1}^2} &\to z_c, 
\end{align}
Since $W_n$ is centred Gaussian $\mathcal{N}_{H_1}(0,T),$ \eqref{eq:existence_quantile} can be expressed as
\begin{equation}
\label{eq:existence_eigenrepresentation}
\frac{r_{n,\gamma}^2 - \mathbb{E}V_n - \sum_{i=1}^{n-1} b_{n,i}^2}
{\sqrt{\operatorname{var} V_n + 4\sum_{i=1}^{n-1} \tau_{i,n}^2b_{n,i}^2}}
\to z_c.
\end{equation}
Here $\{b_{n,i}\}$ has exactly one nonzero entry depending on the smoothness cases $\beta \leq 2\alpha +2p+1$ and $\beta > 2\alpha +2p+1$. The nonzero entry, which we call $b_{n,i_n}$, has the following representation, with $d_n$ to be yet determined,
\begin{align*}
b_{n,i_n}^2 = r_{n,\gamma}^2 -\mathbb{E}V_n - d_n\operatorname{sd}V_n.
\end{align*}
Since $r_{n,\gamma}^2, \mathbb{E}V_n$ and $r_{n,\gamma}^2 - \mathbb{E}V_n$ have the same order and $\operatorname{sd} V_n$ is of strictly smaller order, one can show that the lefthand side of \eqref{eq:existence_eigenrepresentation} is equivalent to
\[
\frac{d_n \operatorname{sd}V_n}{\sqrt{\operatorname{var}V_n + 4 \tau_{i_n,n}^2(r_{n\gamma}^2 - \mathbb{E}V_n)(1+o(1))}},
\]
for bounded or slowly diverging $d_n$. Then \eqref{eq:existence_eigenrepresentation} can be obtained by discussing different smoothness cases separately, by a suitable choice of $i_n,d_n$.

To prove the asymptotic normality in \eqref{eq:existence_normal}, the numerator can be written as
\[
\|W_n + b_n\|_{H_1}^2 - \mathbb{E}\|W_n+b_n\|_{H_1}^2 = \sum_i \tau_{i,n}^2(Z_i^2 - 1) + 2 b_{n,i_n} \tau_{i_n,n}Z_{i_n}.
\]
Next one applies the arguments as in \cite{knapik2011bayesianmild}.

\subsection{Proof of Theorem~\ref{thm:credibleset_extremely}}
\label{sec:proof_credibleset_extreme}

This proof is almost identical to the proof of Theorem 2.2 in \cite{knapik2013bayesianextreme}. We supply the main steps.

Following the same arguments as in the proof of Theorem~\ref{thm:credibleset_mildly}, we obtain 
\begin{align*}
\mathbb{E} X_n &\asymp \rho_n^2(\log (\rho_n^2 n))^{-2\alpha/s} 
\gg 
\operatorname{sd} X_n \asymp \rho_n^2(\log (\rho_n^2 n))^{-1/(2s) - 2\alpha/s} ,\\
\mathbb{E} V_n &\asymp \rho_n^2(\log (\rho_n^2 n))^{-1/s-2\alpha/s} \asymp \operatorname{sd} V_n,
\end{align*}
as in the proof of Theorem 2.2 in \cite{knapik2013bayesianextreme}. This leads to
\[
r_{n,\gamma}^2 \asymp \rho_n^2(\log (\rho_n^2 n))^{-2\alpha/s},
\]
and furthermore,
\[
\mathbb{P}(V_n\leq \delta r_{n,\gamma}^2) 
= 
\mathbb{P}\left(
\frac{V_n - \mathbb{E}V_n}{\operatorname{sd} V_n}
\leq
\frac{\delta r_{n,\gamma}^2 - \mathbb{E}V_n}{\operatorname{sd} V_n}
\right)
\to 1,
\]
for every $\delta>0$.

Similar to Theorem~\ref{thm:credibleset_mildly}, the bounds on the square norm $B_n$ (defined in \eqref{eq:square_norm_of_bias}) of the bias are known: upper bound from the proof of Theorem~\ref{thm:ExtremeIllposed}, and lower bound from Lemma~\ref{lem:supremum_estimation_extreme},
\[
\varepsilon_{n,1}
\lesssim
B_n
\lesssim
\varepsilon_{n,1}\vee\left(\sqrt{n}R \varepsilon_{n,2}\right),
\]
where $\varepsilon_{n,1}, \varepsilon_{n,2}$ are given in \eqref{eq:contraction_rate_extreme}, and $\sqrt{n}R$ satisfies the bound \eqref{eq:remainder_estimation_extreme}.

In case (i.), $B_n\ll r_{n,\gamma}$, and hence \eqref{eq:coverage_B_n_ll_r_n} applies. The rest of the results can be
 obtained in a similar manner.

\section{Auxiliary lemmas}
\label{sec:appendix}
The following lemmas are direct generalisations of the case $s=2$ in the Appendix of \cite{knapik2013bayesianextreme} to a general $s.$ They can be easily proved by simple adjustments of the original proofs in \cite{knapik2013bayesianextreme}, and we only state the results.

\begin{lem}[Lemma 6.1 in \cite{knapik2013bayesianextreme}]
\label{lem:supremum_estimation_extreme}
For $q\in\mathbb{R}$, $u\geq 0$, $v>0$, $t+2q\geq 0$, $p>0$, $0\leq r < pv$ and $s\geq 1$,
\[
\sup_{\|f\|_{S^q}\leq 1} \sum_{i = 1}^{\infty} \frac{f_i^2 i^{-t}e^{-ri^s}}{(1+Ni^{-u}e^{-pi^s})^v}
\asymp
N^{-r/p}(\log N)^{-t/s -2q/s +ru/ps},
\]
as $N\to \infty$.

In addition, for any fixed $f\in S^q$,
\[
N^{r/p}(\log N)^{t/s + 2q/s - ru/ps}\sum_{i = 1}^{\infty} \frac{f_i^2 i^{-t}e^{-ri^s}}{(1+Ni^{-u}e^{-pi^s})^v} \to 0,
\]
as $N\to \infty$.
\end{lem}

\begin{lem}[Lemma 6.2 in \cite{knapik2013bayesianextreme}]
\label{lem:A4}
For $t,u\geq 0$, $v>0$, $p>0$, $0< r < vp$ and $s\geq 1$, as $N\to \infty$,
\[
\sum_{i = 1}^{\infty} \frac{ i^{-t}e^{-ri^s}}{(1+Ni^{-u}e^{-pi^s})^v}
\asymp
N^{-r/p}(\log N)^{-t/s +ru/ps}.
\]
If $r = 0$ and $t>1$, while other assumptions remain unchanged,
\[
\sum_{i = 1}^{\infty} \frac{ i^{-t}e^{-ri^s}}{(1+Ni^{-u}e^{-pi^s})^v}
\asymp
(\log N)^{(-t+1)/s}.
\]
\end{lem}

\begin{lem}[Lemma 6.4 in \cite{knapik2013bayesianextreme}]
\label{lem:A5}
Assume $s\geq 1$.
Let $I_N$ be the solution in $i$ to $N i^{-u}e^{-pi^s} = 1$, for $u \geq 0$ and $p>0$. Then
\[
I_N \sim \left(\frac{1}{p} \log N \right)^{1/s}
\]
\end{lem}

\begin{lem}[Lemma 6.5 in \cite{knapik2013bayesianextreme}]
\label{lem:A6}
Let $s\geq 1$. As $K\to \infty$, we have
\begin{enumerate}[(i.)]
\item for $a>0$ and $b\in \mathbb{R}$,
\[
\int_1^K e^{a x^s} x^b \ dx \sim \frac{1}{a s} e^{a K^s} K^{b - s +1};
\]

\item for $a,b,K>0$,
\[
\int_K^\infty e^{-a x^s} x^{-b} \ dx \leq \frac{1}{a s} e^{-a K^s} K^{-b-s+1}.
\]
\end{enumerate}
\end{lem}

\section*{Acknowledgements}
The research leading to the results in this paper has received funding from the European Research Council under ERC Grant Agreement 320637.

\section*{References}

\bibliography{bib}

\begin{thebibliography}{10}
\expandafter\ifx\csname url\endcsname\relax
  \def\url#1{\texttt{#1}}\fi
\expandafter\ifx\csname urlprefix\endcsname\relax\def\urlprefix{URL }\fi
\expandafter\ifx\csname href\endcsname\relax
  \def\href#1#2{#2} \def\path#1{#1}\fi

\bibitem{alquier2011inverse}
P.~Alquier, E.~Gautier, G.~Stoltz, Inverse Problems and High-Dimensional
  Estimation: Stats in the Ch{\^a}teau Summer School, August 31 -- September 4,
  2009, Lecture Notes in Statistics, Springer, 2011.

\bibitem{bissantz2007inverseregularization}
N.~Bissantz, T.~Hohage, A.~Munk, F.~Ruymgaart, Convergence rates of general
  regularization methods for statistical inverse problems and applications,
  SIAM Journal on Numerical Analysis 45~(6) (2007) 2610--2636.

\bibitem{cavalier2008nonparametricinverseproblems}
L.~Cavalier, Nonparametric statistical inverse problems, Inverse Problems
  24~(3) (2008) 034004.

\bibitem{cavalier2002sharp}
L.~Cavalier, A.~Tsybakov, Sharp adaptation for inverse problems with random
  noise, Probability Theory and Related Fields 123~(3) (2002) 323--354.

\bibitem{cohen2004adaptive_galerkin}
A.~Cohen, M.~Hoffmann, M.~Rei\ss, Adaptive wavelet {G}alerkin methods for
  linear inverse problems, SIAM Journal on Numerical Analysis 42~(4) (2004)
  1479--1501.

\bibitem{donoho1995WVD}
D.~L. Donoho, Nonlinear solution of linear inverse problems by
  wavelet--vaguelette decomposition, Applied and Computational Harmonic
  Analysis 2~(2) (1995) 101--126.

\bibitem{kaipio2006statistical}
J.~Kaipio, E.~Somersalo, Statistical and {C}omputational {I}nverse {P}roblems,
  Applied Mathematical Sciences, Springer New York, 2006.

\bibitem{kirsch2011introduction}
A.~Kirsch, An {I}ntroduction to the {M}athematical {T}heory of {I}nverse
  {P}roblems, Applied Mathematical Sciences, Springer, 2011.

\bibitem{wahba1977integraloperator}
G.~Wahba, Practical approximate solutions to linear operator equations when the
  data are noisy, SIAM Journal on Numerical Analysis 14~(4) (1977) 651--667.

\bibitem{natterer2001mathematics_CT}
F.~Natterer, The {M}athematics of {C}omputerized {T}omography, Classics in
  Applied Mathematics, Society for Industrial and Applied Mathematics, 2001.

\bibitem{isakov2013inverse}
V.~Isakov, Inverse {P}roblems for {P}artial {D}ifferential {E}quations, Applied
  Mathematical Sciences, Springer New York, 2013.

\bibitem{colton2012inverse_acoustic_electromagnetic}
D.~Colton, R.~Kress, Inverse {A}coustic and {E}lectromagnetic {S}cattering
  {T}heory, Applied Mathematical Sciences, Springer New York, 2012.

\bibitem{birke2010inverse_regression}
M.~Birke, N.~Bissantz, H.~Holzmann, Confidence bands for inverse regression
  models, Inverse Problems 26~(11) (2010) 115020.

\bibitem{bissantz2012inverse_regression}
N.~Bissantz, H.~Dette, K.~Proksch, Model checks in inverse regression models
  with convolution-type operators, Scand. J. Stat. 39~(2) (2012) 305--322.

\bibitem{aad:book:2017}
S.~Ghosal, A.~van~der Vaart, Fundamentals of {N}onparametric {B}ayesian
  {I}nference, Vol.~44 of Cambridge Series in Statistical and Probabilistic
  Mathematics, Cambridge University Press, Cambridge, 2017.

\bibitem{knapik2011bayesianmild}
B.~T. Knapik, A.~W. van~der Vaart, J.~H. van Zanten, Bayesian inverse problems
  with {G}aussian priors, Ann. Statist. 39~(5) (2011) 2626--2657.

\bibitem{knapik2013bayesianextreme}
B.~T. Knapik, A.~W. van~der Vaart, J.~H. van Zanten, Bayesian recovery of the
  initial condition for the heat equation, Communications in Statistics --
  Theory and Methods 42~(7) (2013) 1294--1313.

\bibitem{vandervaart2000posterior}
S.~Ghosal, J.~K. Ghosh, A.~W. van~der Vaart, Convergence rates of posterior
  distributions, Ann. Statist. 28~(2) (2000) 500--531.

\bibitem{conway1990courseonFA}
J.~Conway, A {C}ourse in {F}unctional {A}nalysis, Graduate Texts in
  Mathematics, Springer, 1990.

\bibitem{haase2014functional}
M.~Haase, Functional {A}nalysis: {A}n {E}lementary {I}ntroduction, Graduate
  Studies in Mathematics, Amer. Mathematical Society, 2014.

\bibitem{tsybakov2008introduction}
A.~Tsybakov, Introduction to {N}onparametric {E}stimation, Springer Series in
  Statistics, Springer, 2008.

\bibitem{efromovich1998differentiation}
S.~Efromovich, Simultaneous sharp estimation of functions and their
  derivatives, Ann. Statist. 26~(1) (1998) 273--278.

\bibitem{akansu2010GDFTs}
A.~Akansu, H.~Agirman-Tosun, Generalized discrete {F}ourier transform with
  nonlinear phase, IEEE Transactions on Signal Processing 58~(9) (2010)
  4547--4556.

\bibitem{quarteroni2010numerical}
A.~Quarteroni, R.~Sacco, F.~Saleri, Numerical {M}athematics, Texts in Applied
  Mathematics, Springer, 2010.

\end{thebibliography}

\end{document}